\documentclass[12pt, a4paper]{amsart}
\usepackage{amscd,amsthm,amsfonts,amssymb,amsmath,euscript}
\usepackage[dvips]{graphicx}
\usepackage[dvips]{graphics}
\usepackage[matrix,arrow]{xy}
\usepackage{longtable}
\usepackage{url}

\newcounter{statements}
\newtheorem{theorem}[statements]{Theorem}
\newtheorem{conjecture}[statements]{Conjecture}
\newtheorem{problem}[statements]{Problem}
\newtheorem{question}[statements]{Question}

\newtheorem{corollary}[statements]{Corollary}
\newtheorem{proposition}[statements]{Proposition}
\newtheorem{principle}[statements]{Principle}
\newtheorem{fact}[statements]{Fact}
\newtheorem{optimisticpicture}[statements]{Optimistic picture}

\theoremstyle{definition}
\newtheorem{example}[statements]{Example}
\newtheorem{definition}[statements]{Definition}

\theoremstyle{remark}
\newtheorem{remark}[statements]{Remark}

\newcommand{\ZZ}{{\mathbb Z}}

\newcommand{\PP}{{\mathbb P}}
\newcommand{\CC}{{\mathbb C}}

\newcommand{\BB}{{\mathbb B}}
\newcommand{\TT}{{\mathbb T}}
\newcommand{\Aff}{{\mathbb A}}

\newcommand{\Proof}{{\bf Proof}}

\newfont{\smallskob}{cmbx7 scaled\magstep4}
\newfont{\bigskob}{cmbx12 scaled\magstep4}

\newcommand{\itc}[1]{\textup{#1}}

\newcommand{\pic}{\mathrm{Pic}\,}

\newcommand{\Spec}{\mathrm{Spec}\,}

\newcommand{\tit}{Weak Landau--Ginzburg models for smooth Fano threefolds}

\begin{document}

\begin{title}
\tit
\end{title}

\author{Victor Przyjalkowski}

\thanks{The work was partially supported by grants NSF FRG DMS-0854977, NSF DMS-0854977, NSF DMS-0901330, grants FWF P 24572-N25 and FWF P20778,
RFFI grants 11-01-00336-a, 11-01-00185-a, and 12-01-31012, grants MK$-1192.2012.1$,
NSh$-5139.2012.1$, and AG Laboratory GU-HSE, RF government
grant, ag. 11 11.G34.31.0023.}

\address{Steklov Mathematical Institute, 8 Gubkina street, Moscow 119991, Russia} %

\email{victorprz@mi.ras.ru, victorprz@gmail.com}


\maketitle

%\begin{flushright}
%%\begin{minipage}[c]{2.5in}
%To the memory of my teacher V.\,A.\,Iskovskikh.
%%\end{minipage}
%\end{flushright}

\begin{abstract}
We prove that Landau--Ginzburg models for all 17 smooth Fano threefolds with Picard rank 1 can be represented
as Laurent polynomials in 3 variables exhibiting them case by case.
We check that these Landau--Ginzburg models can be compactified to open Calabi--Yau varieties.
In the spirit of L.\,Katzarkov's program we prove that numbers of irreducible components of the central
fibers of compactifications of these pencils are dimensions of intermediate Jacobians of Fano varieties plus 1.
In particular these numbers do not depend on compactifications.
We state most of known methods of finding Landau--Ginzburg models in terms of Laurent polynomials.
We discuss Laurent polynomial representation of Landau--Ginzburg models of Fano varieties and state some problems related to it.
\end{abstract}

%14M25 Toric varieties, Newton polyhedra
%14H10 Families, moduli (algebraic)

 %14Q10 --- Surfaces, hypersurfaces

%14Q15 --- Higher-dimensional varieties
 %14N10 --- Enumerative problems (combinatorial problems)
 %53D45 --- Gromov-Witten invariants, quantum cohomology, Frobenius manifolds
%14N35 --- Gromov-Witten invariants, quantum cohomology
 %14J30 --- $3$-folds
 %14J45 --- Fano varieties
 %14J70 --- Hypersurfaces
%14J81 --- Relationships with physics
 %14D07   Variation of Hodge structures

\bigskip

\section{Introduction}

Mirror Symmetry conjectures relate symplectic properties of a
variety $X$ with algebro-geometric ones for its mirror symmetry pair
--- a variety $Y$ (or one-parametric family of Calabi--Yau varieties $Y\to \Aff^1$)
and vice-versa, relate algebro-geometric properties of $X$ with
symplectic ones of $Y$.
Homological Mirror Symmetry (see~\cite{Ko94}) treats mirror correspondence in terms of derived categories.
It associate two categories with each variety or family.
Given a symplectic form on $X$ (in our considerations it is an anticanonical form $\omega_X$) one can associate
a so called \emph{Fukaya category} $Fuk\,(X)$ with $X$ whose objects are Lagrangian submanifolds with respect to the symplectic form.
The relative version of this category, \emph{a Fukaya--Seidel category} $FS(Y)$ can be associated with $Y$.
Algebraic side of the picture is presented by a derived category of coherent sheaves $D^b(X)$ for $X$ and a derived category
of singularities $D^b_{sing}(Y)$ for $Y$ --- a direct sum of categories over all fibers whose objects are complexes of coherent sheaves modulo perfect complexes.
Homological Mirror Symmetry conjecture for Fano varieties predicts that for any Fano manifold $X$ there exists a so called \emph{Landau--Ginzburg model} $Y\to \Aff^1$ such that their categories are cross-equivalent: $Fuk\,(X)\simeq D^b_{sing}(Y)$ and $D^b(X)\simeq FS(Y)$.

Homological Mirror Symmetry conjecture is very powerful but unfortunately it is very hard to prove it for particular mirror pairs.
So the natural first step is checking coincidence of some invariants of categories discussed above.
A natural invariant of a category $\mathcal C$ is its \emph{Hochschild cohomology} $HH^*(\mathcal C)$. For Fukaya category
Hochschild cohomology is nothing but \emph{Quantum cohomology} (see Subsection~\ref{subsection:quantum cohomology}).
Scaling the symplectic form one can vary Fukaya categories. In this way one obtain a so called \emph{non-commutative Hodge structure}.
In the similar way one can associate a non-commutative Hodge structure with $D^b_{sing}(Y)$
(more precisely, with each singular fiber of $Y\to \Aff^1$). For definitions and constructions of these
structures see~\cite{KKP08}.

The coincidence of non-commutative Hodge structures is called \emph{Mirror Symmetry conjecture of variations of Hodge structures}.
It
enables
one to translate the mirror correspondence for Fano varieties to a
quantitative level.
We formulate it in the following way.
In a lot of cases one can assume that $Y$ is a
complex torus (see Compactification principle~\ref{principle:compactification}). In this case the complex-valued function can be
represented by a Laurent polynomial, which, under some assumptions,  is called \emph{weak
Landau--Ginzburg model} for $X$. Thus it turns out that the problem of
finding a weak Landau--Ginzburg model for $X$ can be reduced to finding
a certain Laurent polynomial. The particular series combinatorially
constructed by this polynomial (the so called \emph{constant terms series})
should be equal to the so called \emph{constant term of regularized $I$-series}
for $X$ constructed by geometrical data (numbers of rational curves
lying on $X$). For more detailed background see~\cite{Prz08}.

There are 17 families of smooth
Fano threefolds with Picard rank $1$, see~\cite{Isk77}. Each of them is determined
by its index and its anticanonical degree.
In this paper we find weak Landau--Ginzburg models (some of them are known but had not been written down) for all 17 families.
It turns out that these models are Laurent polynomials in 3 variables that have \emph{Calabi--Yau compactifications}
to families of K3 surfaces (Theorem~\ref{theorem:main}).

Which numerical invariants of a Fano variety can be reconstructed from
its Landau--Ginzburg model (or from a weak one) and
how? C.\,van Enckevort and D.\,van Straten (\cite{vEvS06}) suggest
to extract characteristic numbers of a general anticanonical section
of Fano variety by writing down a monodromy of the dual family in a
specific basis.

L.\,Katzarkov's recent idea (see say~\cite{KKP08},~\cite{ILP11},~\cite{KP12}) %and~\cite{KP08})
is to
relate the Hodge type of a Fano variety to a structure of a cental
fiber of dual Landau--Ginzburg model and to a sheaf of vanishing
cycles for this fiber.

Theorem~\ref{theorem: number of components} says that numbers
of irreducible components (without multiplicities) of central fibers of {Calabi--Yau
compactifications} of weak Landau--Ginzburg models for Fano threefolds we found
are dimensions of their intermediate Jacobians plus 1.
Actually, under natural assumptions this number does not depend on particular weak Landau--Ginzburg model, see Subsection~\ref{subsection:local to global}
and~\cite{DKLP} (this statement is not clear for Fano threefolds of Picard rank greater than 1, see~\cite{DKLP} for this case).
For rational Fano threefolds this theorem appears in~\cite{AAK12} via constructing Landau--Ginzburg models for these varieties.
For an explanation of phenomenon of Theorem~\ref{theorem: number of components} see Remark~\ref{remark:GKR}.

This paper is the first step in studying weak Landau--Ginzburg models for Fano threefolds.
We refer to~\cite{ILP11},~\cite{DKLP},~\cite{CKP12a},~\cite{CKP12b} where certain properties of some weak Landau--Ginzburg models and their relations to
Homological Mirror Symmetry are studied.

\medskip

We write down a table with weak Landau--Ginzburg models for Fano threefolds here for convenience.
$N$ in the table stays for the number of variety with respect to lexicographic order (index, degree),
$I$ stays for the index of a variety and $\mathrm{deg}$ stays for its anticanonical degree.
Later we prove that the polynomials from the table are actually weak Landau--Ginzburg models
for corresponding Fano varieties (Theorem~\ref{theorem:main}). We also observe most known methods of finding
weak Landau--Ginzburg models (Section~\ref{section:methods of finding}) and discuss some %conjectures and
problems related to them (Section~\ref{section:problems}).
Polynomials in the table are not unique weak Landau--Ginzburg models for Picard rank 1 Fano threefolds (cf. Subsection~\ref{subsection:local to global}).
More examples see, for instance, in~\cite{CKP12a},~\cite{CKP12b}.

%\medskip

\newpage

%\begin{center}
\begin{longtable}{||c|c|c|c|c|p{5cm}|p{5cm}||}
%\caption{Weak Landau--Ginzburg models for Fano threefolds.}
  \hline
  $N$ & $I$ & $\mathrm{deg}$ & $h^{12}$ &
  $\sharp$ comp.  &
\begin{minipage}[c]{5cm}
\centering
\vspace{.1cm}

  Description

\vspace{.1cm}

\end{minipage}
  &
\begin{minipage}[c]{5cm}
\centering
\vspace{.1cm}

  Weak LG model

\vspace{.1cm}

\end{minipage}
\\
  \hline
  \hline
%  \endfirsthead
  1 & 1 & 2 & 52 & 53 &
\begin{minipage}[c]{5cm}
\vspace{.1cm}

Sextic double solid $X_2$ (double cover of $\PP^3$ ramified over smooth sextic).

\vspace{.1cm}

\end{minipage}
&
\begin{minipage}[c]{5cm}
\vspace{.1cm}
\centering

$\frac{(x+y+z+1)^6}{xyz}$

\vspace{.1cm}

\end{minipage}
    \\
  \hline
  2 & 1 & 4 & 30 & 31 &
\begin{minipage}[c]{5cm}
\vspace{.1cm}

  The general element of the family is quartic $X_4$.

\vspace{.1cm}

\end{minipage}
&
\begin{minipage}[c]{5cm}
\vspace{.1cm}
\centering

  $\frac{(x+y+z+1)^4}{xyz}$

\vspace{.1cm}

\end{minipage}
  \\
  \hline
  3 & 1 & 6 & 20 & 21 &
\begin{minipage}[c]{5cm}
\vspace{.1cm}

Smooth complete intersection of quadric and cubic $X_6$.

\vspace{.1cm}

\end{minipage}
&
\begin{minipage}[c]{5cm}
\vspace{.1cm}
\centering

$\frac{(x+1)^2(y+z+1)^3}{xyz}$

\vspace{.1cm}

\end{minipage}
\\
  \hline
  4 & 1 & 8 & 14 & 15 &
\begin{minipage}[c]{5cm}
\vspace{.1cm}

  Smooth complete intersection of three quadrics $X_8$.

\vspace{.1cm}

\end{minipage}
&
\begin{minipage}[c]{5cm}
\vspace{.1cm}
\centering

$\frac{(x+1)^2(y+1)^2(z+1)^2}{xyz}$

\vspace{.1cm}

\end{minipage}
\\
  \hline
  5 & 1 & 10 & 10 & 11 &
\begin{minipage}[c]{5cm}
\vspace{.1cm}

  The general element of the family is $X_{10}$, a section of $G(2,5)$ by 2 hyperplanes and quadric in Pl\"{u}cker embedding.

%\vspace{.1cm}

\end{minipage}

&
%  $\frac{(x^2+x+z+zx+y+yx+yz)^2}{xyz}$ \\
\begin{minipage}[c]{5cm}
\vspace{.1cm}
\centering

$\frac{(1+x+y+z+xy+xz+yz)^2}{xyz}$

%\vspace{.1cm}

\end{minipage}
\\
  \hline
  6 & 1 & 12 & 7 & 8 &
\begin{minipage}[c]{5cm}
\vspace{.1cm}

  Variety $X_{12}$.

\vspace{.1cm}

\end{minipage}
&
\begin{minipage}[c]{5cm}
\vspace{.1cm}
\centering

  $\frac{(x+z+1)(x+y+z+1)(z+1)(y+z)}{xyz}$

\vspace{.1cm}

\end{minipage}
\\
  \hline
  7 & 1 & 14 & 5 & 6 &
\begin{minipage}[c]{5cm}
\vspace{.1cm}

  Variety $X_{14}$, a section of $G(2,6)$ by 5 hyperplanes in Pl\"{u}cker embedding.

\vspace{.1cm}

\end{minipage}
  &

\begin{minipage}[c]{5cm}
\vspace{.1cm}
\centering

  $\frac{(x+y+z+1)^2}{x}$

  $+\frac{(x+y+z+1)(y+z+1)(z+1)^2}{xyz}$

\vspace{.1cm}

\end{minipage}
\\
  \hline
  8 & 1 & 16 & 3 & 4 &
\begin{minipage}[c]{5cm}
\vspace{.1cm}

  Variety $X_{16}$.

  \vspace{.1cm}

\end{minipage}
  &
  \begin{minipage}[c]{5cm}
\vspace{.1cm}
\centering

$\frac{(x+y+z+1)(x+1)(y+1)(z+1)}{xyz}$

\vspace{.1cm}

\end{minipage}
\\
  \hline
  9 & 1 & 18 & 2 & 3 &
\begin{minipage}[c]{5cm}
\vspace{.1cm}

  Variety $X_{18}$.

\vspace{.1cm}

\end{minipage}
  &
\begin{minipage}[c]{5cm}
\vspace{.1cm}
\centering

  $\frac{(x+y+z)(x+xz+xy+xyz+z+y+yz)}{xyz}$

\vspace{.1cm}

\end{minipage}
\\
  \hline
  10 & 1 & 22 & 0 & 1 &
\begin{minipage}[c]{5cm}
\vspace{.1cm}

  Variety $X_{22}$.

\vspace{.1cm}

\end{minipage}
  &
\begin{minipage}[c]{5cm}
\vspace{.1cm}
\centering

  $\frac{xy}{z}+\frac{y}{z}+\frac{x}{z}+x+y+\frac{1}{z}+4$ $+\frac{1}{x}+\frac{1}{y}+z+\frac{1}{xy}+\frac{z}{x}+\frac{z}{y}+
\frac{z}{xy}$

\vspace{.1cm}

\end{minipage}
\\
  \hline
  11 & 2 & $8\cdot 1$ & 21 & 22 &
\begin{minipage}[c]{5cm}
\vspace{.1cm}

  Double Veronese cone $V_{1}$
  (double cover of the cone over the Veronese surface branched in a smooth cubic).

\vspace{.1cm}

\end{minipage}
&
\begin{minipage}[c]{5cm}
\vspace{.1cm}
\centering

%  $\frac{(x^2+y^2+z^2+1)^3}{xyz}$ \\
  $\frac{(x+y+1)^6}{xy^2z}+z$

\vspace{.1cm}

\end{minipage}
\\
  \hline
  12 & 2 & $8\cdot 2$ & 10 & 11 &
  \begin{minipage}[c]{5cm}
\vspace{.1cm}

Quartic double solid $V_2$ (double cover of $\PP^3$ ramified over smooth quartic).

\vspace{.1cm}

\end{minipage}
&
\begin{minipage}[c]{5cm}
\vspace{.1cm}
\centering

$\frac{(x+y+1)^4}{xyz}+z$

\vspace{.1cm}

\end{minipage}
\\
  \hline
  13 & 2 & $8\cdot 3$ & 5 & 6 &
  \begin{minipage}[c]{5cm}
\vspace{.1cm}

Smooth cubic $V_3$.

\vspace{.1cm}

\end{minipage}
&
\begin{minipage}[c]{5cm}
\vspace{.1cm}
\centering

$\frac{(x+y+1)^3}{xyz}+z$

\vspace{.1cm}

\end{minipage}
\\
  \hline
  14 & 2 & $8\cdot 4$ & 2 & 3 &
\begin{minipage}[c]{5cm}
\vspace{.1cm}

  Smooth intersection of two quadrics $V_4$.

\vspace{.1cm}

\end{minipage}
&
\begin{minipage}[c]{5cm}
\vspace{.1cm}
\centering

  $\frac{(x+1)^2(y+1)^2}{xyz}+z$

\vspace{.1cm}

\end{minipage}
\\
  \hline
  15 & 2 & $8\cdot 5$ & 0 & 1 &
\begin{minipage}[c]{5cm}
\vspace{.1cm}

  Variety $V_{5}$, a section of $G(2,5)$ by 3 hyperplanes in Pl\"{u}cker embedding.

\vspace{.1cm}

\end{minipage}
  &
\begin{minipage}[c]{5cm}
\vspace{.1cm}
\centering

  $x+y+z+\frac{1}{x}+\frac{1}{y}+\frac{1}{z}+xyz$

\vspace{.1cm}

\end{minipage}
\\
  \hline
  16 & 3 & $27\cdot 2$ & 0 & 1 &
\begin{minipage}[c]{5cm}
\vspace{.1cm}

  Smooth quadric $Q$.

\vspace{.1cm}

\end{minipage}
&
\begin{minipage}[c]{5cm}
\vspace{.1cm}
\centering

  $\frac{(x+1)^2}{xyz}+y+z$

\vspace{.1cm}

\end{minipage}
\\
  \hline
  17 & 4 & $64$ & 0 & 1 &
\begin{minipage}[c]{5cm}
\vspace{.1cm}

  $\PP^3$.

\vspace{.1cm}

\end{minipage}
&
\begin{minipage}[c]{5cm}
\vspace{.1cm}
\centering

  $x+y+z+\frac{1}{xyz}$

\vspace{.1cm}

\end{minipage}
\\
  \hline
  \hline
\caption[]{\label{table}Weak Landau--Ginzburg models for Fano threefolds.}
\end{longtable}

\section{Mirror Symmetry of variations of Hodge structures}
\label{section:preliminaries}

We consider smooth projective varieties over $\CC$. For any such variety $X$ we denote $H_2(X,\ZZ)/tors$ by $H_2(X)$.
In this paper Calabi--Yau varieties are varieties with trivial canonical class.

\subsection{Regularized quantum $D$-modules}
\label{subsection:quantum cohomology}
Let $X$ be a smooth Fano variety\footnote{This assumption can be weakened; we are interested in the case of
smooth Fano varieties, so we give definitions in this particular case.}. To it one can associate a set of Gromov--Witten invariants
of genus 0. These invariants are numbers counting rational curves lying on $X$.
Consider $\gamma_1,\ldots,\gamma_m\in H^*(X,\ZZ)$, $k_1,\ldots,k_m\in
\ZZ_{\ge 0}$, $m\in \ZZ_+$, and $\beta\in H_2(X)$. The ($m$-pointed) genus 0 Gromov--Witten
invariant with descendants that correspond to this data
(see~\cite{Ma99}, VI--2.1) is denoted by
$$
\langle \tau_{k_1}
\gamma_1,\ldots,\tau_{k_m}\gamma_m\rangle_{\beta}.
$$

Given these invariants (more precisely, prime three-pointed ones, i.e. those with $m=3$ and $k_1=k_2=k_3=0$)
one can define a (small) quantum cohomology ring. This ring is the deformation of the ordinary cohomology ring.

\begin{definition}[see~\cite{Ma99}, Definition $0.0.2$]
Consider a %so called
Novikov ring $\Lambda$ --- the ring of polynomials over $\CC$ in formal variables
$t^\beta$, $\beta\in H_2(X)$, with natural relations $t^{\beta_1}t^{\beta_2}=t^{\beta_1+\beta_2}$.
Fix an effective basis $\Delta\subset H^*(X,\ZZ)$.
The \emph{quantum cohomology ring} is a vector space $QH^*(X)=H^*(X,\ZZ)\otimes \Lambda$ with \emph{quantum multiplication} ---
the bilinear map
$$
\star\colon QH^*(X)\times QH^*(X)\rightarrow QH^*(X)
$$
given by
$$
\gamma_1\star \gamma_2=\sum_{\substack{\gamma\in \Delta,\\ \beta\in H_2(X)}} t^\beta\langle \gamma_1,
\gamma_2, \gamma^\vee\rangle_\beta \gamma
$$
for any $\gamma_1,\gamma_2\in H^*(X)$, where $\gamma^\vee$ is
the Poincar$\mathrm{\acute{e}}$ dual class to $\gamma$ (we identify
elements $\gamma\in H^*(X)$ and $\gamma\otimes 1\in QH^*(X)$).
\end{definition}

Notice that $QH^*(X)$ is graded by $\deg t^\beta=-K_X\cdot \beta$
and the constant term of $\gamma_1\star \gamma_2$ (with respect to
$t$) is $\gamma_1\cdot \gamma_2$, so $QH^*(X)$ indeed is a deformation of $H^*(X, \CC)$.

\medskip

Let $H=-K_X$ and let $QH_H^*(X)$ be the minimal subring of $QH^*(X)$ containing $H$.
Assume it is generated over $\Lambda$ by the linear space $H_H^*(X)$, i.e. $QH^*_H=H_H^*(X)\otimes\Lambda$.
The variety $X$ is called \emph{quantum minimal} if $\dim_{\Lambda} QH_H^*=\dim_\CC H^*_H(X)=\dim X+1$.
Examples of quantum minimal varieties are complete intersections in (weighted) projective spaces
or Fano threefolds with Picard rank 1.

Next we describe the construction of \emph{regularized quantum $\mathcal D$-module} (or,
equivalently, Dubrovin's second structural connection).
More precisely, this $\mathcal D$-module contains an essential submodule corresponding
to $H_H^*(X)$. As we need only this essential part, we give the definition of this submodule;
one should replace $H_H^*(X)$ by $H^*(X)$ in the definition to get the definition of the whole module.
For more particular definition of this $\mathcal D$-module for quantum minimal case see~\cite{GS07},~\cite{Prz07}, and~\cite{Prz08}.

Consider a torus
$\TT=\Spec B$, where $B=\CC[t, t^{-1}]$.
Let $HQ$ be the trivial vector bundle over $\TT$
with fiber $H_H^*(X)$. Let $S=H^0(HQ)$ and let $\star\colon S\times S\to
S$ be the quantum multiplication (we can consider it as an operation on $S\cong QH^*_H(X)\otimes
\CC[t,t^{-1}]$).
%Let $\{h^i\}=\{H^i\otimes 1\}$ be the basis of $S$ over $B$.
Let $\mathcal D=B[\frac{\partial}{\partial t}]$ and
$D=t\frac{\partial}{\partial t}$. Consider a \itc{(}flat\itc{)}
connection $\nabla$ on $HQ$ defined on the sections $\gamma\in H_H^*(X)$ as
$$
\left( \nabla(\gamma),t\frac{\partial}{\partial t}\right)=K_V\star
\gamma
$$
(the pairing is the natural pairing between differential forms and
vector fields). This connection provides the structure of a $\mathcal
D$-module for $S$ by $D(\gamma)=(\nabla(\gamma),D)$.

Let $Q$ be this $\mathcal D$-module. It is not regular in general, so we need to ``regularize'' it to obtain the regular one.
%Let $\gm=\Spec [t,t^{-1}]$.
Let $E=\mathcal D%_\gm
/\mathcal D%_\gm
(t\frac{\partial}{\partial t}-t)$ be the exponential $\mathcal D%_\gm
$-module. %Consider the inclusion $\ZZ(-K_X)\hookrightarrow \pic
%X$. The natural isomorphism $\pic (X)\cong \mathrm{NS}(X)$ ($X$ is
%Fano) and double dualization provide the morphism $j\colon \gm \to
%\TT_{\mathrm{NS}^\vee}$.
Define \emph{the regularization} of $Q$ as
%$Q^{\mathrm{reg}}=j^*(\mu_*(Q\boxtimes j_*(E)))$
$Q^{\mathrm{reg}}=\mu_*(Q\boxtimes E)$, where $\mu\colon
\TT%\gm
\times
\TT%\gm
\to
\TT%\gm
$ is the multiplication and $\boxtimes$ is the
external tensor product. In other words, $Q^{\mathrm{reg}}$ is a convolution
with the anticanonical exponential $\mathcal D$-module. It can be
represented by a differential operator, which is divisible by $D$ on the left: $Q^{\mathrm{reg}}\cong\mathcal D%_{\TT\gm}
/\mathcal
D%_{\gm}
(D
%t\frac{\partial}{\partial t}
L_X)$.
%We denote $t\frac{d}{dt}$ also by $D$.
The differential operator $L_X$ is called \emph{the regularized quantum differential operator}.
If $X$ is quantum minimal, then $L_X$ is said to be~\emph{of type
$DN$}, see~\cite{Go05}, 2.10, and these are studied in~\cite{GS07}. For $N=3$ this operator is given explicitly
in~\cite{Go05}, Example $2.11$, in terms of structural
constants of quantum multiplication by the anticanonical class (two-pointed Gromov--Witten invariants). %If $X$ is a threefold, then
%it is quantum minimal by dimensional reasons, so there corresponds
%an operator of type $D3$ to it. More precisely see in~\cite{Go05}.
Thus, there is an operator of type $D3$ associated to every smooth
Fano threefold with Picard group $\ZZ$. For all smooth Fano threefolds with Picard rank 1 the operators of
type $D3$ are known (see, for instance,~\cite{Go05}, 5.8).

Let $H^0$ be the class in $H^0(X,\ZZ)$ dual to the
fundamental class of a quantum minimal variety $X$. Consider a series
\begin{gather*}
I^X_{H^0}=1+\sum_{\beta} \langle\tau_{-K_X\cdot\beta-2} H^0\rangle_\beta
\cdot t^{-K_X\cdot\beta},
\end{gather*}
where the sum is taken over all $\beta\in H_2(X)$ such that $-K_X\cdot\beta\geq 2$.
For quantum minimal variety $X$ this series is a unique analytic solution
of the equation $L_XI=0$ of type
$$
I^X_{H^0}=1+a_1t+a_2t^2+\ldots\in \CC [[t]],\ \ \ \ a_i\in \CC
$$
(see~\cite{Prz07}, Corollary $2.2.6$, and references therein).

\begin{definition} This series %(The unique) analytic solution of
%$L_XI=0$ of type
%$$
%I^X_{H^0}=1+a_1t+a_2t^2+\ldots\in \CC [[t]],\ \ \ \ a_i\in \CC,
%$$
is called \emph{the fundamental term of the regularized $I$-series}
of $X$.
\end{definition}

\subsection{Weak Landau--Ginzburg models}
Consider a torus $\TT_{LG}=\mathbb
G_{\mathrm{m}}^n=\prod_{i=1}^n \Spec \CC[x_i^{\pm 1}]$ and a
function $f$ on it. This function can be represented by Laurent
polynomial: $f=f(x_1,x_1^{-1}\ldots,x_n,x_n^{-1})$. Let $\phi_f(i)$
be the constant term (i.\,e. the coefficient at $x_1^0\cdot \ldots
\cdot x_n^0$) of $f^i$, and put
$$
\Phi_f=\sum_{i=0}^\infty \phi_f(i)\cdot t^i\in \CC[[t]].
$$

\medskip

\begin{definition} The series $\Phi_f$ %=\sum_{i=0}^\infty \phi_f(i)\cdot t^i$
is
called \emph{the constant terms series} of $f$.
\end{definition}

\begin{definition}
\label{Definition: weak LG models} Let $X$ be a
smooth $n$-dimensional %quantum minimal
Fano variety and
let $I^X_{H^0}\in \CC[[t]]$ be its fundamental term of regularized
$I$-series. The Laurent polynomial $f\in \CC[\ZZ^n]$ is called
\emph{a very weak Landau--Ginzburg model} for $X$ if (up to a shift $f\mapsto f+\alpha$, $\alpha\in \CC$)
$$
\Phi_f=I^X_{H^0}.
$$

A very weak Landau--Ginzburg model $f\in \CC[\ZZ^n]$ is called \emph{a weak
one} if there is a fiberwise compactification of a family $f\colon (\CC^*)^n\to \CC$
whose total space is (an open) smooth Calabi--Yau
variety. Such compactification is called \emph{a Calabi--Yau compactification}.
\end{definition}

\begin{remark}
If the total space of a family and its base are smooth then the general fiber is smooth.
\end{remark}

\begin{remark}
There is a slightly different definition of weak Landau--Ginzburg model in the literature (see say~\cite{Prz08}).
By this definition a weak Landau--Ginzburg model is a very weak one whose general fiber is birational to Calabi--Yau
variety while our definition says that (strengthened) this property holds for all fibers (even reducible ones). However in practice these conditions
are equivalent: we do not know natural examples when these definitions differ.
\end{remark}

The meaning of the definition is the following (see~\cite{BS85},
$10$ or~\cite{Be83}, pp. 50--52). Consider %functions $F_t=1-t\cdot
%f\in \CC[x_1^{\pm 1},\ldots,x_n^{\pm 1}][t]$. They provide
a pencil $\TT_{LG}\to \BB=\PP[u:v]\setminus (0:1)$ with fibers
$Y_\alpha=\{1-\alpha f=0\}$, $\alpha\in \CC\setminus\{0\}\cup \{\infty\}$. %($t$ is the local coordinate around $(1:0)$).

\medskip

The following proposition is a mathematical folklore (see~\cite{Prz08} for the proof).

\begin{proposition}
\label{Proposition: Picard--Fuchs} Assume the Newton polytope of $f\in
\CC[\ZZ^n]$ contains $0$ in its interior. Let $t\in \BB$ be the local
coordinate around $(0:1)$. Then there is a fiberwise $(n-1)$-form
$\omega_t\in \Omega^{n-1}_{\TT_{LG}/\BB}$ and a \itc{(}locally
defined\itc{)} fiberwise $(n-1)$-cycle $\Delta_t$ such that
$$
\Phi_f=\int_{\Delta_t}\omega_t.
$$
\end{proposition}

This means that $\Phi_f$ is a solution of the Picard--Fuchs
equation for the pencil $\{Y_t\}$.

\begin{remark}
\label{remark:equations} Let $PF_f=PF_f(t,
\frac{\partial}{\partial t})$ be a Picard--Fuchs operator of
$\{Y_t\}$. Denote the order of $PF_f$ by $m$ and denote the degree with
respect to $t$ by $r$. Let $Y$ be a semistable compactification of
$\{Y_t\}$ (so we have the map $\widetilde{f}\colon Y\to \PP^1$;
denote it for simplicity by $f$). Denote the dimension of
the transcendental part of $R^{n-1}f_!\,\ZZ_Y$ by $m_f$ (for an algorithm for
computing it see~\cite{DH86}), and denote the number of
singularities of $f$ counted with multiplicities by $r_f$.
Then $m\leq
m_f$ and $r\leq r_f$.
So we can write a differential operator of bounded order by $t$ and $D$ as an operator
with indeterminant coefficients. As $\Phi_f$ annihilates it, we get a system of infinite number of linear equations.
To check that $L_X=PF_f$ we need to solve this system (it has a unique solution, up to scaling, so
we need to solve a finite system of linear equations).

However in practice it is enough to compare the first few coefficients of the
expansion of $\Phi_f$ and $I_{H^0}^X$.
Indeed, we know that the first few terms of $I^X_{H^0}$ determines $L_X$. So if these terms coincide with the first few terms of
$\Phi_f$ then the differential operator vanishing $\Phi_f$ is $L_X$ (up to high order coefficients) which means that $PF_f=L_X$.

\end{remark}

\begin{question}
Let us be given a polytope $\Delta$. Can one find (effectively) a number $s=s(\Delta)$ such that for any polynomial
$f$ whose Newton polytope is $\Delta$, the first $s$ coefficients of $\Phi_f$ determine the other ones
(so in order to prove that $I_{H^0}^X=\Phi_f$ for any $I_{H^0}^X$ it is enough to check coincidence of the first
$s$ coefficients of both series).
In other words, is it true that the linear system of equations on coefficients of a Picard--Fuchs operator
for Laurent polynomial with given Newton polytope for the first $s$ terms is nondegenerate?
\end{question}

\section{Methods for finding weak Landau--Ginzburg models}
\label{section:methods of finding}

We observe here some methods for finding very weak Landau--Ginzburg models for some Fano varieties, in particular for
complete intersections in projective spaces and Grassmannians, and for varieties admitting small toric
degenerations. They are weak ones in practice %--- their general fibers are birational to Calabi--Yau varieties
and usually it is not complicated to prove this. However we do not know a general method of proving this in some cases.

In the next section we find weak Landau--Ginzburg models for Fano threefolds with Picard number 1.
%and discuss their properties.

\subsection{Small toric degenerations}
\label{subsection:small degenerations}
We start from description of mirrors for Fano varieties admitting the so called small toric degenerations.
This description was suggested in~\cite{BCFKS97} and~\cite{Ba97}.

Assume that a smooth Fano variety $X$ admits a degeneration to a terminal Gorenstein toric variety $Y$.
Let $\{v_1,\ldots,v_n\}\subset \ZZ^k$, $v_i=(v_i^1,\ldots,v_i^k)$ be the set of integral generators of rays of the fan of $Y$.
Denote $x^{v_i}=x_1^{v_i^1}\cdot \ldots \cdot x_k^{v_i^k}\in \CC[x_1^{\pm 1},\ldots, x_k^{\pm 1}]$.
Then the weak Landau--Ginzburg model for $X$ (up to a shift $f\to f+\alpha$, $\alpha\in \CC$) is
$$
\sum_{i=1}^n x^{v_i}.
$$

\begin{problem}
\label{problem:small degenerations}
Prove this. In a case of high index one can use Quantum Lefschetz--type arguments (applied to hyperplane section)
and Givental's formula for $I$-series of smooth toric varieties (cf.~\cite{Ga07}, Proposition~1.7.15) for proving
that these polynomials are very weak Landau--Ginzburg models. For proving that %general
%elements of
corresponding families
%are birational to Calabi--Yau varieties
have Calabi--Yau compactifications
one should check that singularities of %these %general
elements of families either
admit a crepant resolution or ``come from the ambient
toric variety''. This is enough for the proof, because these elements are anticanonical sections of ambient toric variety.
For references see~\cite{Ba93}, in particular see Theorem~4.1.9.
\end{problem}

\subsection{Complete intersections}
\label{subsection:complete intersections}
The suggestions for Landau--Ginzburg models for complete intersections in
projective spaces were given by Hori and Vafa in~\cite{HV00}.
In terms of Laurent polynomials their suggestions can be stated in the following way.
The weak Landau--Ginzburg model for smooth complete intersection $X$  of $r$ hypersurfaces
of degrees $k_1,\ldots,k_r$ in $\PP^N$ is (up to the shift $f\to f+\alpha$, $\alpha\in\CC$)
$$
f_X=\frac{\prod_{i=1}^{r}
(x_{i,1}+\ldots+x_{i,k_{i}-1}+1)^{k_i}}{\prod{x_{i,j}\cdot \prod
y_i}}+y_1+\ldots+y_{k_0} \in \CC[\{x_{jl}^{\pm 1},y_s^{\pm 1}\}],
$$
where $k_0=N-\sum k_i$.

\begin{proposition}
\label{proposition:complete intersections}
The polynomial $f_X$ is a weak Landau--Ginzburg model for $X$.
\end{proposition}

\Proof.
According to the well-known Givental's formula for constant term of
$I$-series of $X$ (up to the shift),
$$
I^X_{H^0}=\sum_{i=0}^{\infty} \frac{\prod_{j=0}^r(k_ji)!}{(i!)^{N+1}}t^{k_0i}.
$$
One can check that the constant term of $f_X^n$ is
$\frac{\prod_{j=0}^r(k_ji)!}{(i!)^{N+1}}$ if $n={k_0i}$ and $0$ in the other case.   %\qed

%\begin{remark}
%\label{remark:complete intersections}
Consider the compactification of the pencil corresponding to $f_X$ given by the natural embedding
$(\CC^*)^{N-r}\hookrightarrow \PP(y_0:\ldots:y_{k_0})\times
\PP(x_{1,1}:\ldots:x_{1,k_1})\times\ldots\times\PP(x_{r,1}:\ldots:x_{r,k_r})$.
%The general element
This compactification is (singular) relative hypersurface of multidegree $(k_0+1,k_1,\ldots,k_r)$ in
$\PP^{k_0}\times\PP^{k_1-1}\times\ldots\PP^{k_r-1}$ and hence has trivial canonical class.
%The compactification of the pencil corresponding to $f_X$ given by the natural embedding
%$(\CC^*)^{N-r}\hookrightarrow \PP(y_0:\ldots:y_{k_0})\times
%\PP(x_{1,1}:\ldots:x_{1,k_1})\times\ldots\times\PP(x_{r,1}:\ldots:x_{r,k_r})$ shows that the general element of the pencil
%is birational to (singular) hypersurface of multidegree $(k_0+1,k_1,\ldots,k_r)$ in
%$\PP^{k_0}\times\PP^{k_1-1}\times\ldots\PP^{k_r-1}$ and hence has trivial canonical class.
It is easy to check that singularities of the compactified pencil %general fiber
are ``purely canonical'', which %in all cases we consider,
means that it %elements of these pencils
admits a crepant resolution. Indeed, they are products of du Val singularities and
linear spaces of codimension 2 or products of ordinary double points and linear spaces of codimension 3
(away from intersections of components of singularities). %and they ``intersects transversally'', that is, a
If we blow up any cDV singularity of codimension 2 we get singularities of the same type again
%After blowing up one of intersected components
%another one is du Val singularity along a linear space in the neighborhood of an exceptional divisor
(singularities ``intersect transversally''). Blowing up singularities one by one and taking small resolutions we get a crepant resolution.
Hence the total space of the resolution of the compactification %a general fiber
is a Calabi--Yau variety, so the Laurent polynomial we consider is actually
a weak Landau--Ginzburg model.
\qed

\begin{remark}
\label{remark:weighted projective}
The same can be done for smooth complete intersections of Cartier divisors in weighted projective spaces. %(see~\cite{Prz10}).
That is one can define Hori--Vafa-type Laurent polynomials for them similar to ones for usual intersections complete defined above
and prove that they are very weak Landau--Ginzburg models (\cite{Prz05} and \cite{Prz10}, Theorem 9). For Fano index one case they are known to
be weak ones (see Remark~\ref{remark:X2}). In general case this is not proven yet.
\end{remark}

\subsection{Grassmannians}
Weak Landau--Ginzburg models for Grassmannians were suggested in~\cite{EHX97}, B 25.
Later in~\cite{BCFKS97}, using the construction of small toric degenerations of
Grassmannians (see references therein), these Landau--Ginzburg models were obtained
via small toric degenerations technics.

The suggested model for $G(k,N)$ is
$$
X_{1,1} + \sum_{\substack{1\leq a\leq N-k,\\ 1\leq b \leq k}}
X_{ab}^{-1}(X_{a+1,b} +
X_{a,b+1}) + X_{N-k,k}^{-1}\in \CC[\{X_{ab}, X_{ab}^{-1}\}]
$$
(the variables are $X_{a,b}$, $1\leq a\leq N-k$, $1\leq b \leq k$; for $a>N-k$ or $b>k$ we put $X_{ab}=0$).

\begin{problem}
\label{problem:grassmannians}
Prove that the polynomial written down above is actually weak (or at least very weak) Landau--Ginzburg
model for $G(k,N)$. The problem is combinatorial: to solve it one should find
the coefficients of constant terms series for this polynomial and compare this series with the constant
term of regularized $I$-series for $G(k,N)$ found in~\cite{BCK03} (see also~\cite{BCFKS98}, Conjecture 5.2.3).
For methods for proving the Calabi--Yau condition see Problem~\ref{problem:small degenerations}.
\end{problem}

\begin{remark}
In the similar ways one can write down very weak Landau--Ginzburg models for complete flag manifolds (\cite{Gi96}) and partial flag manifolds (\cite{BCFKS98}). For general flag varieties $G/P$ see~\cite{Re07}. For minuscule varieties see~\cite{BG}.
\end{remark}

\subsection{Complete intersections in Grassmannians}
In this subsection we describe the suggestions for weak Landau--Ginzburg models of
complete intersections in Grassmannians.
The idea for writing them down is the particular case of the method %of writing down of Landau--Ginzburg models
for complete intersections in varieties admitting
small toric degenerations.
More precisely, let $Y=X\cap Y_1\cap\ldots\cap Y_s$ be a complete intersection in a variety admitting a degeneration to terminal Gorenstein
toric variety $X$.
Let $\{v_1,\ldots,v_n\}\subset \ZZ^k$ %, $v_i=(v_i^1,\ldots,v_i^k)$
be the set of integral generators of rays of the fan of $X$ as before.
Let $D(v_i)$'s be divisors corresponding to $v_i$'s.
Let $\{p_1^l,\ldots,p_{r_i}^l\}$, $l=1,\ldots,s$, be subsets of $\{v_1,\ldots,v_n\}$ such that $Y_l=\sum_j D(p_j^l)$ as cohomological classes.
Then the Landau--Ginzburg model for $Y$ is (conjecturally)
$$
\{%\sum_m x^{v_m}=1,\
x^{p_1^l}+\ldots+x^{p_{r_i}^l}=1\}\in \CC[x_1^{\pm 1},\ldots, x_k^{\pm 1}], \ \ l=1,\ldots,s,
$$
with potential $\sum_{i=1}^n x^{v_i}$ (see~\cite{Ba97}).

We describe this procedure for complete intersections in Grassmannian $G(m,r)$ following~\cite{BCFKS97}
(and changing the notation for simplicity)\footnote{It was suggested for Calabi--Yau complete intersections
but it works for Fano varieties in absolutely the same way
(as usual, modulo shift $f\to f+\alpha$, $\alpha\in \CC$).}.
Consider the following Laurent polynomials %$x^{p_1^i}+\ldots+x^{p_{r_i}^i}$, $i=1,\ldots,r$,
in variables $X_{ab}$, $1\leq a \leq r-m$, $1\leq b \leq m$.
$$
  X_{11},
$$
$$
  \frac{X_{1i}}{X_{1,i-1}}+\frac{X_{2i}}{X_{2,i-1}}+\ldots+\frac{X_{r-m,i}}{X_{r-m,i-1}}, \ \ i=2,\ldots,m,
$$
$$
  \frac{X_{j+1,1}}{X_{j1}}+\frac{X_{j+1,2}}{X_{j2}}+\ldots+\frac{X_{j+1,m}}{X_{jm}}, \ \ j=1,\ldots,r-m-1,
$$
$$
\frac{1}{X_{r-m,m}}.
$$
Given any of these polynomials consider the set of rays of the fan of $X$ associated with summands of the polynomial. The sum of the boundary
divisors associated with these rays is equivalent (in cohomology) to a Picard group generator for $X$.
The cohomological class of a hypersurface in Grassmannian is given by its degree.
This means that one can find $s$ sums of boundary divisors such that each divisor is contained in at most one sum and $i$th sum is equivalent to
$Y_i$. This gives us a Landau--Ginzburg model for $Y$.

\begin{problem}
\label{problem:grassmann complete intersections}
In all cases we consider there are (birational) changes of variables for the Landau--Ginzburg models for Fano complete intersections
in $G(m,r)$ obtained in the way described above, such that after these changes Landau--Ginzburg models are
functions on a complex torus (that is, Laurent polynomials).
More precisely, one can express one variable in terms of others from each equation, put them in the Laurent polynomial of Grassmannian
and make a linear change of variables to make denominators monomial. Prove this in the general case.
\end{problem}

It does not follow from this procedure that the Laurent polynomial obtained in this way is a very weak
Landau--Ginzburg model for a complete intersection in $G(m,r)$.
In practice it is even a weak one; one can check this in each particular case.

\begin{problem}
Prove this in the general case (cf. Problems~\ref{problem:small degenerations} and~\ref{problem:grassmannians}).
\end{problem}

\begin{example}
\label{example:V14}
Consider a Fano threefold $X_{14}$. By definition it is the section of $G(2,6)$ cut out by five hyperplanes.
The Landau--Ginzburg model is the variety
$$
\left\{X_{11}=1, \frac{X_{21}}{X_{11}}+\frac{X_{22}}{X_{12}}=1,
\frac{X_{31}}{X_{21}}+\frac{X_{32}}{X_{22}}=1,
\frac{X_{41}}{X_{31}}+\frac{X_{42}}{X_{32}}=1,
\frac{1}{X_{42}}=1
\right\}\subset \CC[\{X_{ij}, X_{ij}^{-1}\}],
$$
$1\leq i\leq 4, \ 1\leq j\leq 2$,
with potential
$$
X_{11}+\frac{X_{21}+X_{12}}{X_{11}}+\frac{X_{22}}{X_{12}}
+\frac{X_{31}+X_{22}}{X_{21}}+\frac{X_{32}}{X_{22}}
+\frac{X_{41}+X_{32}}{X_{31}}+\frac{X_{42}}{X_{32}}
+\frac{X_{42}}{X_{41}}+\frac{1}{X_{42}}.
$$

Denote $X_{12}=a$, $X_{22}=b$, $X_{32}=c$.
Then
$$
X_{21}=\frac{a-b}{a},
$$
$$
X_{31}=\frac{(a-b)(b-c)}{ab},
$$
$$
X_{41}=\frac{(c-1)(a-b)(b-c)}{abc}.
$$
So the potential is
$$
5+a+\frac{ab}{a-b}+\frac{abc}{(a-b)(b-c)}+\frac{abc}{(a-b)(b-c)(c-1)}=5+\frac{a^2}{a-b}+\frac{abc^2}{(a-b)(b-c)(c-1)}.
$$
Denote $x=a-b$, $y=b-c$, $z=c-1$. Then $a=x+y+z+1$, $b=y+z+1$, $c=z+1$ and we get the potential
$$
f=5+\frac{(x+y+z+1)^2}{x}+\frac{(x+y+z+1)(y+z+1)(z+1)^2}{xyz}\in \CC[x,x^{-1},y,y^{-1}, z, z^{-1}].
$$
The constant term of the regularized $I$-series for $X_{14}$ (shifted by 4) is
$$
I_{H^0}^{X_{14}}=1+4t+48t^2+760t^3+13840t^4+273504t^5+5703096t^6+\ldots
$$
(see~\cite{Prz04}). It is easy to see that the constant terms series for $f_{14}=f-5$ equals
$I_{H^0}^{X_{14}}$ up to more then $16$ coefficients. This means that they are equal (see Remark~\ref{remark:equations}).

The fiberwise compactification of %the general element of
the pencil $f=\lambda$, $\lambda\in \CC$ given by the natural map
$\Spec \CC[x,x^{-1},y,y^{-1}, z, z^{-1}]\to \PP(x:y:z:w)$ gives a family of quartic surfaces. Singularities of this
family are ordinary double points or du Val
along lines (see the proof of Proposition~\ref{proposition:complete intersections}).
So the family admit a Calabi--Yau compactification.

Therefore the Laurent polynomial we obtain is a weak Landau--Ginzburg model for $X_{14}$.
\end{example}

\section{Fano threefolds of Picard rank 1}
\label{section:Fano threefolds}
In this section we study Laurent polynomials from Table~\ref{table} and prove that they are weak Landau--Ginzburg models
for corresponding Fano varieties.
First we describe how these weak Landau--Ginzburg models (we call them \emph{standard}\footnote{This term is local for the
paper. Weak Landau--Ginzburg models from Table~\ref{table}, from our point of view, are not better then others.}) are obtained.

\medskip

Varieties $\bf 1$, $\bf 11$, $\bf 12$
are hypersurfaces in weighted projective spaces of degree 6 in $\PP(1:1:1:1:3)$, of degree 6 in $\PP(1:1:1:2:3)$,
and of degree 4 in $\PP(1:1:1:1:2)$ respectively. Their weak Landau--Ginzburg
models can be found by Hori--Vafa procedure similar to procedure for complete intersections in projective
spaces described in Subsection~\ref{subsection:complete intersections} (see Remarks~\ref{remark:X2}
and~\ref{remark:weighted projective}). %and Remark~\ref{remark:veronese}).

Varieties $\bf 2$, $\bf 3$, $\bf 4$, $\bf 13$, $\bf 14$, $\bf 16$, $\bf 17$ are complete intersections, so
their weak Landau--Ginzburg models can be found using Proposition~\ref{proposition:complete intersections}.

Varieties $\bf 5$, $\bf 7$, $\bf 15$ are complete intersections in Grassmannians, so the corresponding polynomials can be obtained using
the method described in Problem~\ref{problem:grassmann complete intersections}.
The polynomial for $X_{14}$ is studied in Example~\ref{example:V14}.
There is another way to obtain the same polynomials for $V_5$ and $X_{10}$.
Indeed, $V_5$ has a small toric degeneration (this is proved by S.\,Galkin in his Thesis~\cite{Ga07}), so its weak Landau--Ginzburg
model is given by the corresponding polytope (see~\ref{subsection:small degenerations}).
According to V.\,Golyshev (see~\cite{Go05}), %and~\cite{ILP11},
the Landau--Ginzburg model for $X_{10}$ is a quotient of the model for $V_5$
(see~\cite{ILP11} and~\cite{DKLP} for the proof).
Taking invariants of the quotient and changing coordinates one can get a weak Landau--Ginzburg model for $X_{10}$;
the form we write down is convenient for calculations.

The polynomial for $\bf 6$-th variety $X_{12}$ is found in~\cite{BP84}.
We change coordinates a bit to get the convenient form as written.

Finally, polynomials for varieties $\bf 8$, $\bf 9$, $\bf 10$ are found in~\cite{Prz08}.
There is a misprint in the polynomial for $X_{16}$ in the journal version of~\cite{Prz08}; it is  corrected in the preprint on arXiv.
It is remarkable that some of these polynomials were found under the assumption that there are %canonical
Gorenstein toric degenerations of corresponding varieties. Later S.\,Galkin in his Thesis (\cite{Ga07}) proved that there is a terminal Gorenstein
toric degeneration of $X_{22}$, so the corresponding polynomial can be obtained using a method from Subsection~\ref{subsection:small degenerations}.

\begin{theorem}
\label{theorem:main}
Standard polynomials are weak Landau--Ginzburg models for Fano threefolds with Picard rank 1.
\end{theorem}

\Proof. Direct computations show that these polynomials are very weak Landau--Ginzburg models (see Remark~\ref{remark:equations}).
Straightforward compactifications  of pencils corresponding to all polynomials except for polynomials for $V_1$ and $X_{2}$
give relative quartics in $\PP^3$ or (for complete intersections) relative anticanonical sections in products
of projective spaces. %Remark~\ref{remark:complete intersections})
Singularities of obtained total spaces %of obtained threefolds
are du Val singularities
along lines and ordinary double points on each step of their minimal resolutions
(cf. Proposition~\ref{proposition:complete intersections}). This gives Calabi--Yau compactifications of families.
%they are birational to K3 surfaces.

%Let us prove that %a general fiber
%fibers of the
Consider a
pencil
$$
\frac{(x+y+1)^6}{xy^2z}+z=\lambda,\ \ \ \ \lambda\in \CC,
$$
associated with $V_1$. %can be compactified to %is birational to a
%K3 surfaces.
Compactify this family to a family of surfaces in $\Aff^3$:
$$
(x+y+1)^6=(\lambda-z)xy^2z.
$$
Then changing the variables $a=x+y+1$, we get
$$
a^6=(\lambda-z)(a-y-1)y^2z.
$$
Consider this family as a family of surfaces lying in $(\CC^*)^3$ with coordinates $a$, $y$, $z$ (in other words let us divide by $a$, $y$, $z$).
Changing the variables $b=y/a$, $c=yz/a$, we get a family
$$
a^4=(\lambda b-c)(a-ab-1)c.
$$
Compactify it to a family in $\PP^3$. We get a family of quartics with du Val singularities along lines
(see the proof of Proposition~\ref{proposition:complete intersections}). %, whose resolution is a K3 surface.
After a resolution of these singularities one get a family whose total space is Calabi--Yau and
the initial family of hypersurfaces in torus is embedded to this resolution.
%Thus the total space of its resolution is a Calabi--Yau variety.

Finally %prove that %a general element
%elements of the
consider a
pencil for $X_2$. %are birational to a K3 surface.
We have a family
$$
(x+y+z+1)^6=\lambda xyz.
$$
In birational coordinates $a=x+y+z+1$, $b=x/(x+y+z+1)$, $c=y/(x+y+z+1)$ we get
$$
a^4=\lambda bc(a-ab-ac-1).
$$
As before one can check that the initial family can be compactified to a Calabi--Yau threefold.
\qed

\begin{remark}
\label{remark:X2}
Let us prove that the general element of a standard weak Landau--Ginzburg model for $X_2$ can be compactified to
a K3 surface in another, more conceptual way (suggested by V.\,Golyshev).
Remember that a Hori--Vafa mirror for a hypersurface of degree $d$ in $\PP(w_0:\ldots:w_n)$
is
$$
\left\{%
\begin{array}{ll}
    y_0^{w_0}\cdot\ldots\cdot y_n^{w_n}=1\\
    y_0+\ldots+y_k=1,\\
\end{array}%
\right.
$$
%$$
%y_0^{w_0}\cdot\ldots\cdot y_n^{w^n}=1,
%$$
%$$
%y_0+\ldots+y_k=1,
%$$
where $w_0+\ldots+w_k=d$,
with the potential
$$
f=y_0+\ldots+y_n.
$$
For $X_2$, the hypersurface of degree $6$ in $\PP(1:1:1:1:3)$,
we have (up to a shift $f\to f-1$) the variety
$$
\left\{%
\begin{array}{ll}
y_0y_1y_2y_3y_4^3=1\\
y_1+y_2+y_3+y_4=1
\end{array}%
\right.
$$
with the potential
$$
f=y_0.
$$
Taking change of variables
$$
y_1=\frac{x}{x+y+z+t}, \ \ y_2=\frac{y}{x+y+z+t},\ \ y_3=\frac{z}{x+y+z+t},\ \ y_4=\frac{t}{x+y+z+t}
$$
(where $x,y,z,t$ are projective coordinates)
we get the Landau--Ginzburg model
$$
y_0xyzt^3=(x+y+z+t)^6, \ \ \ \ f=y_0.
$$
So, in local chart, say $t\neq 0$, we finally get the weak Landau--Ginzburg model
$$
f_2^\prime=\frac{(x+y+z+1)^6}{xyz}.
$$
The general element of the pencil corresponding to $f_2^\prime$ is birational to the general element of the initial
Hori--Vafa model. Inverse the potential: $u=1/f$. Then we get the pencil
$$
y_1y_2y_3y_4^3=u,\ \ \ \ y_1+y_2+y_3+y_4=1.
$$
This model is exactly the Landau--Ginzburg model for $\PP(1:1:1:3)$ (see~\cite{CG06}, (2)). So, by
Theorem 1.15 in~\cite{CG06}, a general element of the pencil we are interested in is birational to a K3 surface.
\end{remark}

\section{Reconstructing Hodge numbers}

The following fact can be proved immediately by constructing Hironaka's house.

\begin{fact}
Let $X$ and $Y$ are two birational smooth Calabi--Yau varieties. Then they are birational in codimension 1.
\end{fact}

By definition of a flop they differ by flops.

\begin{corollary}
\label{corollary:Calabi--Yau}
Any two Calabi--Yau compactifications of weak Landau--Ginzburg model differ by flops.
\end{corollary}

\begin{theorem}
\label{theorem: number of components} Let $X$ be a smooth Fano
threefold with Picard rank 1. Let $f\colon (\CC^*)^3\to\CC$ be its
standard weak Landau--Ginzburg model. Let $k_X$ be the number of irreducible
components of the central fiber of a Calabi--Yau compactification of $f$.
Then $k_X=h^{12}(X)+1$.
\end{theorem}

\Proof.
The proof is
given by direct calculations of Calabi--Yau compactifications considered in Theorem~\ref{theorem:main}
in all 17 cases, one by one (see Examples~\ref{example:V16} and~\ref{example:V18} and~\cite{IKP11},~\cite{KP12}).
By Corollary~\ref{corollary:Calabi--Yau}
%Theorem~\ref{theorem: minimal compactification}
numbers of components of central fibers
of all Calabi--Yau compactifications of $f$ are the same.
\qed

\begin{remark}
As $X$ is a Fano threefold and $\pic (X)=\ZZ$, the Hodge diamond of $X$ is the following.
\begin{center}
\begin{tabular}{ccccccc}
   &   &     & 1 &     &   &   \\
   &   &  0  &   &  0  &   &   \\
   & 0 &     & 1 &     & 0 &   \\
0  &   & $h^{12}(X)$ &   & $h^{12}(X)$ &   & 0 \\
   & 0 &     & 1 &     & 0 &   \\
   &   &  0  &   &  0  &   &   \\
   &   &     & 1 &     &   &   \\
\end{tabular}
\end{center}
Thus Theorem~\ref{theorem: number of components} enables one to reconstruct all Hodge numbers of Picard rank 1 Fano threefolds.
\end{remark}

\begin{remark}
\label{remark:GKR}
According to~\cite{GKR12}, the phenomenon appeared in Theorem~\ref{theorem: number of components} can be explained (at least for complete
intersections) as follows. Consider a perverse sheaf $\mathcal F$ of vanishing cycles to the central fiber of $f$.
If Homological Mirror Symmetry holds for its fiberwise compactification, then $k_X$ can be computed via spectral sequence from~\cite{GKR12} for
$H^1(\mathcal F)$ and $H^3(\mathcal F)$.
Hence by~\cite{GKR12} $H^1(\mathcal F)$ and $H^3(\mathcal F)$ are isomorphic as Hodge structures to $H^{12}(X)$ and $H^{21}(X)$.
From this point of view the statement of Theorem~\ref{theorem: number of components} holds for any Landau--Ginzburg model
(cf. Corollary~\ref{corollary:uniqueness}). %, even of any dimension.
Thus Theorem~\ref{theorem: number of components} confirms Homological Mirror Symmetry conjecture for Picard rank 1 Fano
threefolds and their (compactified) standard Landau--Ginzburg models.
\end{remark}

One can see on examples that this phenomenon holds in higher Picard rank case.

\begin{question}
Direct computations, even in Picard rank 1 case as in Theorem~\ref{theorem: number of components}, can be very complicated (usually for big Hodge numbers).
So doing them for all Picard ranks is a huge technical problem. Is it possible to do it (or at least prove Theorem~\ref{theorem: number of components})
in more conceptual way?
\end{question}

\begin{remark}
An analog of Theorem~\ref{theorem: number of components} in higher dimensions, due to~\cite{GKR12},
is the following. Let us have an $n$-dimensional Fano variety $X$. Let $k$ be a number of components
of a central fiber of its (toric) Landau--Ginzburg model of the same dimension. Then
Homological Mirror Symmetry involves that $k=h^{1,n-1}(X)+1$.
\end{remark}

\begin{problem}
Prove this for toric Landau--Ginzburg models of Fano varieties of dimension greater then 3.
\end{problem}

\begin{example}
Consider a general cubic fourfold. Its Hodge structure is a sum of a part generated by a hyperplane section and a primitive Hodge
structure in dimension 4 --- weight 2 Hodge structure of dimensions $(1,20,1)$. In particular $h^{13}=1$ and $h^{22}=21$.
The central fiber of Calabi--Yau compactification of Hori--Vafa type toric Landau--Ginzburg model has 2 components intersecting
by particular Picard rank 20 K3 surface (see~\cite{KP09}).
%Is it true that for $n$-dimensional Fano variety a number of central components is $h^{1,n-1}$?
\end{example}

\begin{example}
\label{example:V16}
Consider the variety $X_{16}$. %(8-th in Table~\ref{table}).
Its standard weak
Landau--Ginzburg model is
$$
\frac{(x+1)(y+1)(z+1)(x+y+1)}{xyz}.
$$
The compactification in a projective space is a family of quartics
$$
\{(x+t)(y+t)(z+t)(x+y+z+t)=\lambda xyzt\}\subset \Aff[\lambda]\times
\PP[x:y:z:t].
$$
There are 4 components of the central fiber $\lambda=0$.
Singularities are the disjoint union of 9 ``horizontal'' lines. All
of them are products of du Val singularities of type $A_1$ and affine lines in the
neighborhood of a general point. After blowing them up we get three
ordinary double points in the central fiber. So finally we get no
new components of the central fiber and $k_{X_{16}}=4=h^{12}(X_{16})+1$.
\end{example}

\begin{example}
\label{example:V18}
Consider the variety $X_{18}$. % (9-th in Table~\ref{table}).
Its standard weak
Landau--Ginzburg model is
$$
\frac{(x+y+z)(x+xz+xy+xyz+z+y+yz)}{xyz}.
$$
The compactification in a projective space is a family of quartics
$$
\{(x+y+z)(xt^2+xzt+xyt+xyz+zt^3+yt^2+yzt)=\lambda xyzt\}\subset
\Aff[\lambda]\times \PP[x:y:z:t].
$$
There are 2 components of the central fiber $\lambda=0$.
Singularities are 3 ``horizontal'' lines globally of type $A_1$
along lines, 3 ``horizontal'' lines globally of type $A_2$ along lines,
and one ``horizontal'' line $\ell$ which is, away from the central fiber, of type $A_1$ along a line.
The intersection of two components of
the central fiber is a plane cubic with one node; this node lies on
$\ell$. Blowing $\ell$ up we get one more ``vertical'' line of singularities globally
of type $A_1$ along a line. Two components of the fiber over 0 intersect
now at the union of two lines (one of them is a singularity of our
threefold); these lines intersect at two points. Blowing the remaining
singularity up we get three surfaces over 0. Each two of them intersect by
a rational curve, and three such lines intersect at two points. So finally we
get $k_{X_{18}}=3=h^{12}(X_{18})+1$.
\end{example}

\begin{remark}
More complicated cases of compactifications are described in
detail in~\cite{IKP11},~\cite{KP12}.
\end{remark}

\section{Properties, examples, problems}
\label{section:problems}

In the previous sections we considered weak Landau--Ginzburg models for Picard rank 1 Fano threefolds.
However for given variety such weak Landau--Ginzburg model is not unique.
In this section we discuss how to choose ``correct'' ones and what depends on a particular choice.
We consider two approaches: ``global to local'' one claiming which weak Landau--Ginzburg models are correct
and ``local to global'' one saying that for Picard rank one Fano threefolds nothing depend on a particular choice. The
second approach is studied and discussed in~\cite{DKLP}.

\subsection{Global to local}
\label{subsection:global to local}
Let us come back to the initial definition of dual Landau--Ginzburg model for a Fano variety $X$ of dimension $n$.
Remind that a Laurent polynomial $f$ in $n$ variables is called a very weak Landau--Ginzburg model for $X$ if $I_{H^0}^X=\Phi_f$.
In other words, \emph{period condition} for $f$ is satisfied.
However this condition is not enough to ``feel geometry'' of $X$. To get stronger definition one can use the following principle.

\begin{principle}[Compactification principle]
\label{principle:compactification}
There exists a fiberwise compactification of the family of fibers of ``good'' very weak Landau--Ginzburg model (defined up to flops) satisfying (B-side of) Homological Mirror Symmetry conjecture.
\end{principle}

In particular this means that it should exist a fiberwise compactification to an (open) smooth Calabi--Yau variety --- family of compact
Calabi--Yau varieties.
This condition is strong enough: say if $f(x_1,\ldots, x_n)$ is a very weak Landau--Ginzburg model for $X$ then for big enough $k$
Laurent polynomial $f(x_1^k,\ldots, x_n)$ is a very weak Landau--Ginzburg model for $X$ but compactification principle fails for it.
Thus a necessary condition for very weak Landau--Ginzburg model to satisfy the compactification principle is \emph{a Calabi--Yau condition}
saying that there should exist a fiberwise compactification whose total space is a smooth Calabi--Yau variety.

However these two conditions are not enough.

\begin{example}
\label{example:4cubic}
A 4-dimensional cubic (see~\cite{KP09}) has the following two weak Landau--Ginzburg models:
$$
\frac{(x+y+1)^3}{xyzw}+z+w
$$
and
$$
\left(x_1+x_2+\frac{1}{x_1x_2}\right)\left(y_1+y_2+\frac{1}{y_1y_2}\right).
$$
The Calabi--Yau compactification of the second one has complicated central fiber with many components and wrong invariants.
%The Hori--Vafa-type polynomials look ``more natural''.
\end{example}

\begin{example}
\label{example:non-minkowski}
Let $X$ be a complete intersection of two quadrics in $\PP^5$. Consider the following weak
Landau--Ginzburg model for $X$:
$$
\left(x+\frac{1}{x}\right)\left(y+\frac{1}{y}\right)\left(z+\frac{1}{z}\right).
$$
The number of components over zero of its compactification in $\PP^3\times \Aff^1$ is 30 while $h^{12}(X)=2$.
\end{example}

\begin{remark}
Coordinates of singular fibers of weak Landau--Ginzburg model are determined by its
Picard--Fuchs equation. Example~\ref{example:non-minkowski} shows that even
for Picard rank 1 Fano threefolds
the number of components of fiber over zero is not determined by an equation.
\end{remark}

The last condition we want to put on Laurent polynomials is related to toric degenerations (cf. Subsection~\ref{subsection:small degenerations}).
Consider a Fano variety $X$ and its weak Landau--Ginzburg model $f$. We say that $f$ satisfies \emph{toric condition} if
there is an embedded degeneration $X\rightsquigarrow T$ to a toric variety $T$ whose fan polytope (the convex hull of integral generators of rays
 of $T$'s fan) coincides with the Newton polytope (the convex hull of non-zero coefficients) of $f$. In this case $f$ is called
\emph{a toric Landau--Ginzburg model}.

One can easily check that weak Landau--Ginzburg models from Examples~\ref{example:4cubic} and~\ref{example:non-minkowski} are not toric:
degrees of Fano varieties differ from degrees of toric varieties associated with Laurent polynomials.

From this point of view we state the following conjecture.

\begin{conjecture}[Strong version of Mirror Symmetry of variations of Hodge structures conjecture]
Any smooth Fano variety has a toric Landau--Ginzburg model.
\end{conjecture}

\begin{corollary}
Any smooth Fano variety has a toric degeneration.
\end{corollary}

By Theorem~\ref{theorem:main}, Proposition~\ref{proposition:complete intersections}, and~\cite{ILP11} this conjecture holds for Picard rank 1 Fano threefolds
and complete intersections.

So we hope that the following picture holds.

\begin{optimisticpicture}
Toric degenerations of smooth Fano varieties are in 1-to-1 correspondence with toric Landau--Ginzburg models.
Compactification principle holds for them.
\end{optimisticpicture}

\begin{question}
Is the opposite to the second part of optimistic picture holds? That is it true that all Landau--Ginzburg models
(from Homological Mirror Symmetry point of view) of the same dimension as an initial Fano variety are compactifications
of toric ones? In particular, is it true that all of them are rational?
\end{question}

These questions are treated in~\cite{DKLP}.

\begin{question}
Should we put some conditions on toric varieties in the optimistic picture?
\end{question}

\begin{remark}
\label{remark:multitoric}
Toric varieties can have several smoothings to different smooth Fano varieties, see~\cite{fanosearch},~\cite{CKP12a}, and~\cite{CKP12b}.
\end{remark}

\subsection{Local to global}
\label{subsection:local to global}
Another way of determining which weak Landau--Ginzburg models are ``correct'' is known only for threefolds. We sketch it here for Picard rank 1 Fano varieties.
For more details see~\cite{DKLP}.

As we mentioned, Homological Mirror Symmetry conjecture expects that fibers of Landau--Ginzburg model for Fano variety are Calabi--Yau varieties.
More precisely, they are expected to be mirror dual to anticanonical sections of a Fano variety. For the threefold case
this duality is essentially the classical Dolgachev--Nikulin duality of K3 surfaces.

Let $H$ be a hyperbolic lattice, $\ZZ\oplus \ZZ$ with intersection form
$$
\left(
  \begin{array}{cc}
    0 & 1 \\
    1 & 0 \\
  \end{array}
\right).
$$
The intersection lattice on the second cohomology on any K3 surface is $K=H\oplus H\oplus H\oplus E_8(-1)\oplus E_8(-1)$.
Consider a family $U_L$ of K3 surfaces whose lattice of algebraic cycles contains $L\subset K$ (and coincides with $L$ for general K3 surface
in the family).
Consider a lattice $M'=L^\bot$, the orthogonal to $L$ in $K$. Let $M'=H\oplus M$.

\begin{definition}
The family of K3 surfaces $U_M$ is called \emph{the Dolgachev--Nikulin dual} family to $U_L$.
\end{definition}

Consider a principally polarized family of anticanonical sections of a Fano threefold $X$ of index $i$ and degree $(-K_X)^3=i^3k$.
It is nothing but $U_{L_n}$ with $L_n=\langle 2n \rangle$, $2n=ik$. The lattice $L_n$ is a sublattice of $H$. Using this embedding to one of the $H$-summands
of $K$ we can see that its Dolgachev--Nikulin dual lattice is $M_n=H\oplus E_8(-1)\oplus E_8(-1)+\langle -2n \rangle$.

The surfaces with Picard lattices $M_n$ are \emph{Shioda--Inose surfaces}. They are resolutions of quotients of specific K3 surfaces $S$
by \emph{Nikulin involution}, the one keeping transcendental lattice $T_S$; it interchanges two copies of $E_8(-1)$. Another description of Shioda--Inose surfaces is Kummer ones
going back to products of elliptic curves with $n$-isogenic ones. $M_n$-polarized Shioda--Inose surfaces form an $1$-dimensional irreducible
family.

It turns out that fibers of standard toric Landau--Ginzburg models can be compactified to Shioda--Inose surfaces dual to anticanonical sections
of Fano threefolds:

\begin{theorem}[\cite{DKLP}]
Let $X$ be a Picard rank 1 Fano threefold of index $i$ and $(-K_X)^3=i^3k$. Then fibers of Calabi--Yau compactifications of standard toric Landau--Ginzburg model for $X$
are Shioda--Inose surfaces with Picard lattice $M_{ik/2}$.
\end{theorem}

We call (toric) weak Landau--Ginzburg models with Shioda--Inose condition \emph{good}.

Thus compactifications of good Landau--Ginzburg models are, modulo coverings and the standard action of $PSL(2,\CC)$ on the base,
the unique families of corresponding Shioda--Inose surfaces. More precisely, they are index-to-one coverings of the moduli spaces.

\begin{corollary}[cf. Corollary~\ref{corollary:Calabi--Yau}]
\label{corollary:uniqueness}
A compactification (to a smooth Calabi--Yau variety) of good weak Landau--Ginzburg model is unique up to flops.
\end{corollary}

This means that if Homological Mirror Symmetry for Picard rank 1 Fano threefolds holds then Landau--Ginzburg models for it are, up to flops,
compactifications of the standard ones. Moreover, all other good Landau--Ginzburg models are birational over the base $\Aff^1$ to the standard
Landau--Ginzburg models.

\subsection{Toric degenerations and Laurent polynomials}
\label{subsection:coefficients}

There are several ways how to find (toric) weak Landau--Ginzburg models for given Fano variety $X$. Some of them
are described in Section~\ref{section:methods of finding}. However in some cases these methods do not work,
say, if an easy geometric description for $X$ is not known. In this case one can hope to ``guess'' weak Landau--Ginzburg
model via toric degenerations in the following way.
First, one should guess a toric variety with the same numerical invariants as $X$ --- degree, Picard rank, etc.
Then (this is not necessary) one should prove that this toric variety is a degeneration of $X$.
Thus, by optimistic picture, there should be a toric Landau--Ginzburg model whose Newton polytope
is a fan polytope of the toric variety we find. The problem is to reconstruct a Laurent polynomial with
prescribed constant terms series from its Newton polytope (``put coefficients on integral points of the polytope'').
This problem is not solved yet. However there are several approaches for doing it, mostly in canonical
and up to threefold case.

\emph{Binomial principle} says that on vertices of polytope we should put 1's and
on the $i$'th (from any end) integral point of an edge of integral length $n$ is $\binom{n}{i}$.
This principle works for a lot of cases (in other words, for toric varieties with cDV singularities
that is, ones whose integral points of a fan polytope (except for the origin) lie on edges).

\begin{conjecture}[Prokhorov,~\cite{Pr05}]
Fano threefolds with cDV singularities are smoothable.
\end{conjecture}

However this is not always the case, and we need to consider worse singularities. For them
(for canonical\footnote{Canonicity means that the origin is a unique integral point in the interior of the fan polytope.} threefold case)
\emph{Minkowski principle}, suggested in~\cite{fanosearch}, can be applied.

This principle is a generalization of binomial principle in the following way. Remind that a segment of integral length $n$
is a Minkowski sum of $n$ indecomposable integral polytopes --- segments of length 1. A unique polynomial in one variable $x$, having
coefficient 1 at vertices of this segment is $x+1$. Binomial principle says that coefficients at length $n$ segment are given
by a product of $n$ polynomials $x+1$ corresponding to each Minkowski summand.
This description can be generalized to two-dimensional case, that is, to facets of 3-dimensional polytope.
Consider each facet and its Minkowski decomposition to indecomposable summands. Assume that all of these summands are of type
$A_k$, $k\geq 0$ --- Newton polytopes of polynomials $(x+1)^k+y$; if this is not the case the principle can't be applied.
Then coefficients for the initial facet are given by product of the polynomials corresponding to summands.

The idea of the Minkowski decomposition of facets is that such decompositions vanish some obstructions for deformations of toric varieties (see~\cite{Alt94}).
This means, in particular, that different Minkowski decompositions of facets give different Laurent polynomials corresponding to different
smooth Fano varieties having degenerations to the same toric variety (see Remark~\ref{remark:multitoric}).
The other corollary is that toric varieties for which Minkowski principle can be applied tend to be smoothable in the anticanonical embedding.
%Say, they can be smoothed to all 98 non-hyperelliptic Fano threefolds.
From the other hand non-Minkowski Laurent polynomials
tend to ``correspond'' (after defining mirror correspondence for them) to toric degenerations in non-anticanonical embedding
or even non-smoothable varieties.
%All canonical very weak Landau--Ginzburg models for smooth Fano threefolds satisfy Minkowski principle (see~\cite{fanosearch}).
%However this principe does not necessary holds if we extend the notion of weak Landau--Ginzburg models to say %orbifolds
%(see for instance~\cite{CG06}).

\begin{example}
Toric Landau--Ginzburg model (from a point of view of~\cite{CG06} or~\cite{AKO04}) for $\PP(1: a_1:\ldots: a_n)$ is
$x_1+\ldots+x_n+\frac{1}{x_1^{a_1}\cdot\ldots\cdot x_n^{a_n}}$.
\end{example}

\begin{example}[\cite{CKP12a},~\cite{CKP12b}]
Consider $\PP(1:1:2:4)$. Non-Minkowski toric Landau--Ginzburg model for $\PP^3$, a Laurent polynomial
$$\frac{(x+1)^2}{xyz}+\frac{y}{z}+z,$$
corresponds to smoothing $\PP(1:1:2:4)$ to $\PP^3$ as quadrics in $\PP(1:1:1:1:2)$
(see~\cite{KP12},
Example~2.13). Another non-Minkowski polynomial, $x+y+z+1/xy^2z^4$, is a Landau--Ginzburg model %(in some sense)
for $\PP(1:1:2:4)$ itself.
\end{example}

\begin{question}
Is it true that for toric Landau--Ginzburg models for smooth Fano varieties associated with toric degenerations in anticanonical embedding
Minkowski principle holds? What about ``if and only if''?
\end{question}

\begin{remark}
Unfortunately, it is not enough to consider Gorenstein\footnote{
Toric variety is Gorenstein if a polytope dual to its fan polytope is integral. In particular, Gorenstein toric varieties are canonical.} toric
Landau--Ginzburg models for all Fano threefolds. For example sextic double solid has no Gorenstein toric Landau--Ginzburg model as
there is no integral polytope of volume $\frac{2}{3!}=\frac{1}{3}$ containing only one integral point in the interior.
For similar reasons there is no Gorenstein toric Landau--Ginzburg model
for double Veronese cone.
Standard toric Landau--Ginzburg models for both of these varieties are not canonical.
All the remaining Picard rank 1 Fano threefolds have Gorenstein toric Landau--Ginzburg models. % by Remark~\ref{remark:1}.
%It seems to be correct to consider polytopes that correspond to
%canonical $\QQ$-Gorenstein toric varieties.
Among 105 smooth Fano threefolds 98 ones can have Gorenstein toric Landau--Ginzburg models. All of them have at least very weak ones
(see~\cite{fanosearch}). For more details and examples see~\cite{CKP12a}.
\end{remark}

\begin{example}
There are 5 possible Gorenstein toric degenerations for $\PP^3$ (see~\cite{Ka08}). Three of them give toric Landau--Ginzburg models.
There are 5 possible Gorenstein toric degenerations for quadric in $\PP^4$ (see~\cite{Ka08}). Four of them give toric Landau--Ginzburg models.
\end{example}

\begin{remark}
There is an infinite number of (non-Gorenstein) toric weak Landau--Ginzburg models for $\PP^3$. Say, ones whose Newton
polytopes are fan polytopes of $\PP(a^2,b^2,c^2,abc)$, where $(a,b,c)$'s are Markov triples, that is, triples
of natural numbers whose greatest common divisor is 1 and $a^2+b^2+c^2=3abc$ (cf.~\cite{HP05} and~\cite{GU}).
Due to Prokhorov's remark $\PP^3$ does degenerate to such projective spaces (this follows from~\cite{HP05}, Theorem 1.1).
\end{remark}

\medskip

The author is grateful to I.\,Cheltsov, S.\,Galkin, V.\,Golyshev, A.\,Iliev, L.\,Katzarkov, V.\,Lazic, V.\,Nikulin,
D.\,Orlov, K.\,Shramov, and A.\,Wilson
for helpful comments and %to S.\,Galkin for
important remarks, and to referee for advises on reorganization of the paper.

\end{document}